\documentclass{article}

\usepackage{amsmath}%
\usepackage{amsfonts}%
\usepackage{amssymb}%
\usepackage{amsthm}
\usepackage{graphicx}
\usepackage{tikz-cd}
\usepackage{hyperref}
\usepackage[T1]{fontenc}

\DeclareUnicodeCharacter{3000}{ }

\theoremstyle{plain}

\newtheorem*{acknowledgement*}{Acknowledgement}

\newtheorem{claim}{Claim}

\newtheorem{lemma-definition}{Lemma-Definition}

\newtheorem{definition}{Definition}

\numberwithin{equation}{section}
\newtheorem{theorem}{Theorem}[section]
\newtheorem{lemma}[theorem]{Lemma}
\newtheorem{proposition}[theorem]{Proposition}
\newtheorem{corollary}[theorem]{Corollary}

\theoremstyle{remark}
\newtheorem{remark}[theorem]{Remark}
\newtheorem{example}[theorem]{Example}
\numberwithin{equation}{section}

\newcommand{\sff}{\mathrm{I\!I}} % second fund form 
\DeclareMathOperator{\tr}{Tr}

\DeclareMathOperator{\prv}{Pr_{ver}}
\DeclareMathOperator{\diag}{diag}

\DeclareMathOperator{\grad}{grad}
\newcommand{\Mn}{\mathcal{M}(n,\mathbb{R})}
\newcommand{\pdef}{\mathcal{S}_+(n,n)}
\DeclareMathOperator{\strat}{\textit{\dj}}

\newcommand{\psdef}{\mathcal{S}_+(k,n)} 
\DeclareMathOperator{\Ad}{Ad}

\title{A model of invariant control system using mean curvature drift from Brownian motion under submersions}
\author{Huang Ching-Peng}

\begin{document}
	\maketitle

\begin{abstract}

	Given a Riemannian submersion $\phi: M \to N$, we construct a stochastic process $X$ on $M$ such that the image $Y:=\phi(X)$ is a (reversed, scaled) mean curvature flow of the fibers of the submersion. The model example is the mapping $\pi: GL(n) \to GL(n)/O(n)$, whose image is equivalent to the space of $n$-by-$n$ positive definite matrices, $\pdef$, and the said flow has deterministic image. We are able to compute explicitly the mean curvature (and hence the drift term) of the fibers w.r.t. this map, (i) under diagonalization and (ii) in matrix entries, writing mean curvature as the gradient of log volume of orbits. As a consequence, we are able to write down Brownian motions explicitly on several common homogeneous spaces, such as Poincar\'e's upper half plane and the Bures-Wasserstein geometry on $\pdef$, on which we can see the eigenvalue processes of Brownian motion reminiscent of Dyson's Brownian motion.
	
	By choosing the background metric via natural $GL(n)$ action, we arrive at an invariant control system on the $GL(n)$-homogenous space $GL(n)/O(n)$. We investigate the feasibility of developing stochastic algorithms using the mean curvature flow.
	
	\vfill
	
\end{abstract}

\tableofcontents

\begin{acknowledgement*}
	Partial support for this work was provided by the National Science Foundation under grant DMS-2107205
\end{acknowledgement*}

\section{Motivation and overview}

All notations and abbreviations are listed in \ref{notations}. Readers can go to the section for quick reference.

We are motivated by the question: how does one write down Brownian motion on the space of positive semidefinite matrices $\psdef$ under the Bures-Wasserstein geometry? This geometry is defined by the Riemannian submersion \begin{equation}\label{BW def 1}
	\begin{array}{cccc}
		\pi_k: & \mathbb{R}^{nk^\times} &  \to  &  \mathcal{S}_+(k,n) \cong \mathbb{R}^{nk^\times}/O(k) \\
		& M & \mapsto  & MM^*
	\end{array},
\end{equation}  where $n\geq k$ and $\mathbb{R}^{nk^\times}$ is the set of all full rank $n\times k$ real matrices as an open submanifold of the Euclidean space $\mathbb{R}^{nk}$.  $\mathcal{S}_+(k,n)$ is the set of all $n$-by-$n$ positive semidefinite matrices of rank $k$. 

On a manifold, one can write down Brownian motion in explicit local coordinates. However in the case of Bures-Wasserstein, such is not a simple task.

Moreover, a local coordinate system is not always the most useful way to describe Brownian motion--practically, for manifolds embedded in Euclidean spaces, it is often favourable to use the embedding global coordinates. For quotient manifolds, the total spaces often have simpler structures, we naturally turn our attention to make relations between the Brownian motion on the total space and that in the quotient. 

\subsection{Brownian motion through submersions and on homogeneous spaces}

	It\^o differential differs from classical calculus in the second derivative term when applying the chain rule. In the language of Riemannian geometry, the second derivative often appears as the mean curvature of a submanifold.

	An important example is 
	\begin{theorem}\label{BMsubmersion}\cite{watson1973manifold}
	Let $\phi: (M,g) \to (N,h)$ be a Riemannian submersion, and $W_M$ is a Brownian motion on $(M,g)$, then the image of $W_M$ under $\phi$ can be written as \begin{equation}
		d\phi(W_M) = -c\phi_*(H_{\phi^{-1}\circ\phi(W_M)  }(W_M ) )dt + dW_N.
	\end{equation}
	Here $c = \frac{\dim(\phi^{-1}\phi(W_M))}{2}$ is a constant depending on the relative dimension of the submersion.  $W_N$ is a Brownian motion on $(N,h)$. The vector $H_{\phi^{-1}\circ\phi(W_M)  }(W_M )$ is the mean curvature of the fibre $\phi^{-1}(M)$ as an embedded submanifold of $M$ at $W_M$. 
\end{theorem} 

The proof is by a straight forward computation of projecting the Laplace-Beltrami operator on the total space to the quotient.  We shall prove a version of the theorem that extracts the drift term. (See \ref{ito vertical BM}.)

We narrow down our focus on the Riemannian submersions coming from actions of Lie groups as isometries. These are of the form
\begin{equation}
	\phi: M \to M/K;
\end{equation} which arise from an action of the Lie group $K$ on $M$ isometrically.  Since the mean curvature $\vec{H}$ only depends on the local metric information, the mean curvature flow (generally a PDE) for a group orbit $K\cdot p$ within $M$  descends to an ODE on $M/K$,
\begin{equation}\label{first time J shows up}
	\frac{\partial \phi(x)}{\partial t} = J(\phi(x)),
\end{equation} where $J := -c\phi_*(\vec{H})$ with the same constant $c$ as in the theorem.

 We explore by examples. First we consider Riemannian homogeneous spaces and in particular symmetric spaces. These cases have totally geodesic fibres, and BM ``commutes'' with the submersion.

The B-W geometry is an example with a nonzero drift term, with $K = O(k)$ as in \ref{BW def 1}. We give specific formulae and include the calculation in \ref{BWBM section}. 

\begin{proposition}\label{mean curvature calculation}
	Let $M\in GL(n)$ and $P = MM^*\in\pdef$.  
	\begin{enumerate}
		\item If $P = U\Lambda U^*$, for some $U\in O(n)$ and $\Lambda = \diag(\lambda_1, \cdots, \lambda_n)$. Then $J := -c\pi_{k*}(\vec{H})$, with $c = \frac{\dim(O(n))}{2}$, is 
		\begin{equation}
			J(P) = U\diag(\sum_{j \neq 1}\frac{\lambda_1}{\lambda_1+\lambda_j}, \cdots, \sum_{j \neq k}\frac{\lambda_k}{\lambda_k+\lambda_j}, \cdots)U^*.
		\end{equation}
		
		\item Alternatively, the same is computed in terms of the entries of $M$ (or $P$) via first writing down the metric tensor $g|_{MO(n)}$ of $MO(n) \subset GL(n)$. 
		\begin{equation}
			\vec{H}(M) = -\nabla \log\det g|_{MO(n)},
		\end{equation}   where an entry of $g|_{MO(n)}$ is (using Kronecker $\delta$) \begin{equation}
			g_{ij,kl} = \delta_{kl}P_{jl} + \delta_{jl}P_{ik} - \delta_{il}P_{jk} - \delta_{jk}P_{il}.
		\end{equation} The gradient $\nabla$ here is taken under the Euclidean metric on $GL(n)$, and  \begin{equation}
			J(P) = -\frac{2}{\dim O(n)} (HM^* + MH^*).
		\end{equation}
	\end{enumerate}
\end{proposition}

 As a corollary, we have the It\^o equation for BM on $(\pdef, BW)$:
 
 \[
 dW_{BW} =  (dW_{GL(n)})W_{GL(n)}^* + W_{GL(n)}dW_{GL(n)}^* +(nI -J(W_{BW}) )dt.
 \]

\subsection{A control system using the mean curvature drift}
	
	Inspired by the mean curvature drift in \ref{mean curvature calculation}, we turn our attention to developing a stochastic algorithm on the B-W geometry using this drift. We start with the nonsingular case $\pdef = GL(n)/O(n)$. The core idea is that we consider the set of metrics on $GL(n)$ defined by 

	\[\left< V,W\right>_{R} := \tr_{R} VW^*: = \tr RVW^*, \] for $R \in \pdef$. 
	
	Secondly, define $\prv$ to be the orthogonal projection to the vertical direction with respect to the submersion and let $B$ be a BM on the total space, we prove that the process defined by the It\^o equation
	\begin{equation}
		dM = \prv dB 
	\end{equation} exists. (See \ref{ito vertical BM}.) Moreover, let $P=\phi(M)$ be its image, \begin{equation}
	dP = J(P)dt,
	\end{equation} with $J= -c\phi_*H$ as before. Roughly speaking, the vertical part of a BM contributes exactly the mean curvature drift. 
	
	Now, choosing a path $R_t$, we  have a corresponding evolution 
	\begin{equation}
		dP = J^{R_t}(P)dt, 
	\end{equation} where $J^R$ is the drift term under the new metric $\tr_{R}$. Thinking $R_t$ as the control, we have a control system that is stochastic on the total space $GL(n)$ but deterministic on the quotient $\pdef$.
	
	The following gives detailed descriptions of the system:

\begin{proposition}\label{control equivalence}
	The following are equivalent descriptions of the same control system for the evolution $P_t \in \pdef \cong GL(n)/O(n)$ in the sense that the controls give the same possible set of $\{P_t\}$ for a fixed time $t$ and an initial condition $P_0$.
	
	\begin{enumerate}
		\item $\dot{P_t} = M_tC_tM_t^*$, where $M_tM_t^* = P$, and the control $C_t$ is any matrix in $\pdef$ such that the  eigenvalues are of the form 
		\begin{equation}
			(\sum_{j \neq 1}\frac{1}{\lambda_1+\lambda_j}, \cdots, \sum_{j \neq k}\frac{1}{\lambda_k+\lambda_j}, \cdots)
		\end{equation} with each $\lambda_j >0$. Here the choice of  $M_t \in M_t O(n)$ does not change the set of possible $\dot{P_t}$. 
		
		\item $\dot{P_t} = J^{R_t}(P_t)$, where the control is $R_t \in \pdef$, and $J^R$ denotes the drift term as previously under metric $\tr_R$.
		
		\item $\dot{P_t} = G_t^{-*}J(G_t^*P_tG_t)G_t^{-1}$, where the control is any $G_t \in GL(n)$.
	\end{enumerate}
\end{proposition}

In particular, (3) in the proposition means what we have here is an \textit{invariant control system}; that is, the accessible directions at $P$, defined as $\mathfrak{a}(P) = \{J^R(P)| R\in \pdef\} \subset T_P \pdef$ is the restriction to $T_P \pdef$ of a set of vector fields $\mathfrak{a}$ on $\pdef$ invariant under the $GL(n)$ action descending from the left $GL(n)$-multiplication on itself. 

Invariant control systems already exist in literature. We estimate the \textit{reachable set}, i.e. the set of points on $\pdef$ the evolution can reach, as follows.

\begin{proposition}\label{main thm} For the control system described in \ref{control equivalence}
	\begin{enumerate}
		
		\item  \begin{equation}
			P_{t_2} - P_{t_1},
		\end{equation} is positive definite if $t_2>t_1$.
		
		\item $\mathfrak{a}(I)$ is open in $T_I\pdef$. Therefore the reachable set via piece-wise constant controls $\mathcal{R}^{pc}(I)$ has nonempty interior.
		
		\item  The reachable set $\mathcal{R}^{reg}(I)$ via \textit{regulated controls} contains a subset 
		\begin{equation}
			I + 	\{U\diag(\exp\lambda_1,\cdots, \exp\lambda_n)U^*| U \in O(n), \; (\lambda_1,\cdot, \lambda_n)\in cone\{e_i + e_j|i\neq j \} \},
		\end{equation} where $\{e_j\}$ denotes the standard basis of $\mathbb{R}^n$. 
	\end{enumerate}
\end{proposition}

	\subsection{Notations and abbreviations}\label{notations}
		\begin{itemize}
		\item Matrices, including matrix groups, are over real numbers if not specified. 
		
		\item The transpose of a real matrix $A$ is denoted as $A^*$. Although throughout the paper, mostly we only consider real matrices, when extending to complex numbers, $AA^*$ is positive definite but not $AA^T$. $(A^{-1})^*$ and $(A^{-1})^T$ are abbreviated as $A^{-*}$ and $A^{-T}$ since the two operations commute. 
		
		\item Manifolds (especially Lie groups) are assumed to be finite dimensional if not specified. 
		
		\item $\pdef$ denotes the space of $n\times n$ positive definite matrices of full rank.   $\mathcal{S}(n)$ is the space of all $n \times n$ (real) symmetric matrices, whereas $\mathcal{S} _+(k,n)$ is the locus of semidefinite matrices of rank $k$. 
		
		\item The differential of a smooth map, $\phi$, between smooth manifolds is written as $\phi_*$ to avoid confusion with stochastic differentials. 
		
		\item For stochastic processes, we often suppress the time variable $t$ for simplicity. For example, we write $X$ instead of $X_t$ when it is clear in the context. 
		
		\item For It\^o differentials we use $d$ while for Stratonovich differential we use $\strat$.
		
		\item\textit{BM} is short for \textit{Brownian motion}. \textit{MCF} stands for \textit{mean curvature flow}, sometimes scaled or reversed if clear from the context.  
		
		\item $\sff(\cdot, \cdot)$ denotes the second fundamental form, and we take the convention that the mean curvature of an embedded manifold $M$ is defined as $\frac{1}{\dim M}\tr \sff$.
		
		\item $\vec{H}(x)$ (sometimes without the vector arrow when not emphasizing it is a vector) denotes the mean curvature vector at a point $x$. In most cases through out the paper, $\vec{H}$ is a vector field. $J(x)$ denotes the image of $\vec{H}(x)$ through a submersion. 
		
		When we need to specific the \textit{metric} under which $\vec{H}$ is derived, we write $\vec{H}^{metric}$ and $J^{metric}$.

		\item  \[I_{k,n} = \begin{bmatrix}
			I_k & \\
			& 0
		\end{bmatrix},\] where $I_k$ is the $k\times k$ identity matrix as usual.
	\end{itemize}

\section{Preliminaries}
	In this section, we organise known results and definitions needed.
	
	\subsection{SDE on manifolds}
	
	\begin{definition}\label{SDE man} (C.f. \cite{hsu2002stochastic} definition 1.2.3)
		
		Let $M$ be a smooth manifold, $\{V_\alpha\}_\alpha$ a set of vector fields on $M$, and $\{Z^\alpha\}_\alpha$ real valued semimartingales. Adopting the Einstein summing convention, an $M$-valued process $X_t$ is a solution to the SDE \[dX = V_\alpha\strat Z^\alpha\] if for any smooth function $f\in \mathcal{C}^\infty(M)$, \[df(X) = V_\alpha f(X)\strat Z^\alpha\] up to some stopping time $\tau$. 
	\end{definition} 

	In particular, we may choose $f$ to be the embedding coordinates or local coordinates. 
	
	Often we are interested in the It\^o form of the equations and need to convert.

	\subsection{Mean curvature flow of group orbits}
	
	There are many examples of Riemannian submersions in general. We narrow down to a class of particular interest. That is, the submersion of the form \begin{equation}
		\phi: M \to N = M/K,
	\end{equation} where $K$ is a Lie group acting on $M$ as a subgroup of the isometries of $M$, then $\phi_*(H)$ is constant along the fibre, i.e. $\phi_*(H(x))  = \phi_*(H(y))\in T_{\phi(x) }N$ for any $y \in \phi^{-1}\circ\phi(x)$.  In this case, the MCF ``commutes'' with the submersion. That is, the mean curvature flow $x_t$, in general cases described by PDEs, of the fibres in $M$ descends to a simple ODE of $p_t:= \phi(x_t)$ on $N$,  \begin{equation}
		\dot{p_t} = -c\phi_*(H(x_t)) =: J(p_t),
	\end{equation} determined solely by the vector field $J$ on $M/K$.

	A particular fact  will be useful for calculating the mean curvature.
	
	\begin{proposition} \label{log vol grad}(C.f. \cite{pacini2003mean} Proposition 1)
		Let $M$ be a Riemannian manifold such that a compact Lie group $K$ acts as a subgroup of isometries of $M$. The volume of orbit $vol(K\cdot m)$ is a smooth function on $M$, and the vector field of orbit mean curvatures is \begin{equation}
			\vec{H}(m) = -\grad\log vol(K\cdot m)
		\end{equation} for $m\in M$.
	\end{proposition}
	
	The volume function above can be obtained from the following recipe:
	
	\begin{lemma}\label{volume of orbit}(C.f. \cite{pacini2003mean} Proof of Proposition 1)
		Let $(M,g)$ be a Riemannian manifold such that a compact Lie group $K$ acts as a subgroup of isometries of $M$, and let $\{Z_\alpha\}_\alpha$ be a basis of $\mathfrak{k} = Lie(K)$, the volume of an orbit is \begin{equation}
			vol(K\cdot m) = \sqrt{\det i^*g}\int_K \wedge_\alpha Z_\alpha^*.
		\end{equation} Here $i^*g$ is the metric, in terms of the basis $\{Z_\alpha\}_\alpha$, pulled back from the inclusion $i: K\cdot m \hookrightarrow M$, and $Z_\alpha^*$ is the dual of $Z_\alpha$. 
	\end{lemma}

	\subsection{Invariant Brownian motion on Lie groups}\label{invariant BM}
	
	\begin{definition}(C.f. \cite{article}) 
		Let $G$ be a finite-dimensional Lie group and $\mathfrak{g}$ its Lie algebra. If $\mathfrak{g}$ is equipped with an inner product, and therefore it is a Euclidean space, then a \textit{left invariant Brownian motion} $X$ on $G$ is defined by
		\[dX = X\strat W := XY_\alpha \strat W^\alpha : =  Y_\alpha|_X\strat W^\alpha ,\] where $Y_\alpha$'s are an orthonormal basis of $\mathfrak{g}$, identified as the vector space of invariant vector fields on $G$, and $W^\alpha$'s are standard real Wiener processes so $W := W^\alpha Y_\alpha$ is a standard (Euclidean) Brownian motion on $\mathfrak{g}$. 
	\end{definition}

	The rightmost form in the definition is the one compliant to Definition \ref{SDE man}. The middle form however helps to consider the equation in Lie-theoretic settings especially when $G$ is a matrix group, as in practice we often identify the Lie algebra as the tangent space at the identity of $G$.
	
	For a \textit{unimodular} Lie group $G$, i.e. left invariant Haar measures on $G$ are right invariant, or equivalently $\det Ad \equiv 1$ on $G$ or $\tr ad \equiv 0$ on $\mathfrak{g}$, the Laplace-Beltrami operator is a ``sum of squares''. In particular, the general linear group $GL(n)$, special linear group $SL(n)$, and compact groups such as $O(n)$ are unimodular and the following apply. 
	
	\begin{theorem}(\cite{urakawa1979least} \S 2)
		Let $G$ be a finite-dimensional Lie group and $\mathfrak{g}$ its Lie algebra equipped with an inner product. If $G$ is unimodular and $\{Y_i\}_i$ is an orthonormal basis of $\mathfrak{g}$, then the Laplace-Beltrami operator under the left invariant metric on $G$ induced from the inner product on $\mathfrak{g}$ is  \[\Delta_G = \sum_i Y_i^2. \]
	\end{theorem}  
	
	\begin{corollary}
		Let $G$ be a finite-dimensional unimodular Lie group and $\mathfrak{g}$ its Lie algebra equipped with an inner product. Then left invariant Brownian motion on $G$ is exactly the Brownian motion under the left invariant metric induced from the inner product of $\mathfrak{g}$.
	\end{corollary}

	\begin{proof}
		Let $\{Y_i\}$ be an orthonormal basis of $\mathfrak{g}$. The $G$-valued process defined by $dW = Y_i\strat W^i$ has infinitesimal generator $\frac{1}{2}\Delta_G = \frac{1}{2}\sum_i Y_i^2$. 
	\end{proof}

\section{Brownian motion on homogeneous spaces}	

With the background knowledge, we may explore some examples. 

	\subsection{It\^o vertical Brownian motion}
	
	We start with a version of Theorem \ref{BMsubmersion}.

	\begin{lemma}\label{ito vertical BM}
		Let $\phi: M \to N$ be a submersion and $M$ be equipped with a metric $g$. Let $\prv$ denote the orthogonal projection to the vertical direction with respect to the submersion. Define the process $X$ by 
		\begin{equation}\label{X process}
			dX = \prv dW, 
		\end{equation} where $W$ is a Brownian motion on $(M,g)$.  Assuming the existence of solutions, the image process, $\phi(X)　=: P \in N$, satisfies the equation \begin{equation}\label{P process}
			dP = -c\phi_*(H(X_t)) dt, 
		\end{equation}where $H(X)$ is the mean curvature of the fiber $\phi^{-1}(P)$ at $X$ in $M$ and $c = \frac{\dim \phi^{-1}(P)}{2}$ depends on the dimension of the fibre.
	\end{lemma}
	
	\begin{remark}
		In the original theorem \ref{BMsubmersion}, we require the submersion $M\to N$ to be Riemannian. For the lemma \ref{ito vertical BM} here, we do not need the metric structure on $N$. 
	\end{remark}
	
	\begin{proof}
		Pick a normal coordinate system at $o := X_0$ such that the first $m$ coordinate derivatives $\partial_1, \cdots, \partial_m$ span the vertical space at $o$.  
		
		Let $W_M$ be a Brownian motion on $M$. In local coordinates (c.f. \cite{hsu2002stochastic} Example 3.3.5), 
		\begin{equation}
			dW_M^k = (\sqrt{g^{-1}}dB)^k - \frac{1}{2}g^{ij}\Gamma_{ij}^k dt.
		\end{equation}  % check
		
		Let $\xi_i = \prv \partial_i$. 
		\begin{equation}
			dX = \xi_i dW_M^i = \xi_i\strat W_M^i - \frac{1}{2}D_{\xi_i} \xi_j d\langle W_M^i, W_M^j \rangle
		\end{equation} in Stratonovich form, where $D$ denotes the ``Euclidean connection'' for the coordinates. 
		
		At $o$, $D_{\xi_i} \xi_j$ coincides with $\nabla_{\partial_i}\partial_j$ and $\langle W_M^i, W_M^j \rangle = t$. So 
		\begin{equation}
			dP = d\phi(X) = - \frac{1}{2}\sum_{i = 1,\cdots, m} \phi_*\sff(\partial_i, \partial_i)dt.
		\end{equation}
		
	\end{proof}

		The use of the lemma will be apparent in Section 4. However we introduce it here for a cleaner illustration in the following example.

	\begin{example}\label{circle example} % add interpretation and the other pic
		Consider the punctured Euclidean space $M = \mathbb{R}^n\setminus 0$ quotient by $K = SO(n)$. The orbits are $n$-spheres. The quotient space is $(0,\infty)$.
		
		\begin{equation}
			dX = \prv(X)dW = (I - \frac{XX^*}{X^*X}) dW.
		\end{equation}
		
		Let $R = \sqrt{X^*X}$ and $S = R^2 = X^*X$. By It\^o calculus, 
		\begin{equation}
			\begin{split}
				dS & = dX^*X + X^*dX + \frac{1}{2}dX^*dX \\
				& = 2X^*dX + \frac{1}{2}dW^*\prv(X)^*\prv(X)dW \\
				& = 0 + \frac{1}{2}dW^*(I - \frac{XX^*}{S})dW \\
				& = \frac{1}{2}(ndt - \frac{1}{S}Sdt )  = \frac{n-1}{2}dt.
			\end{split}
		\end{equation}
		So \begin{equation}\label{radial process squared}
			S = R^2 = \frac{n-1}{2}t + S_0.
		\end{equation}
		
		On the other hand, the mean curvature and $R$ are related as 
		\begin{equation}
			\dot{R} = -\frac{n-1}{2}\vec{H} = \frac{n-1}{2R},
		\end{equation} which leads to the same as \autoref{radial process squared}. We have shown in this example, the It\^o vertical projection of Brownian motion is indeed the mean curvature drift predicted by Theorem \ref{BMsubmersion}.
		
		\begin{figure}[h]
			\centering
			\includegraphics[width=0.7\linewidth]{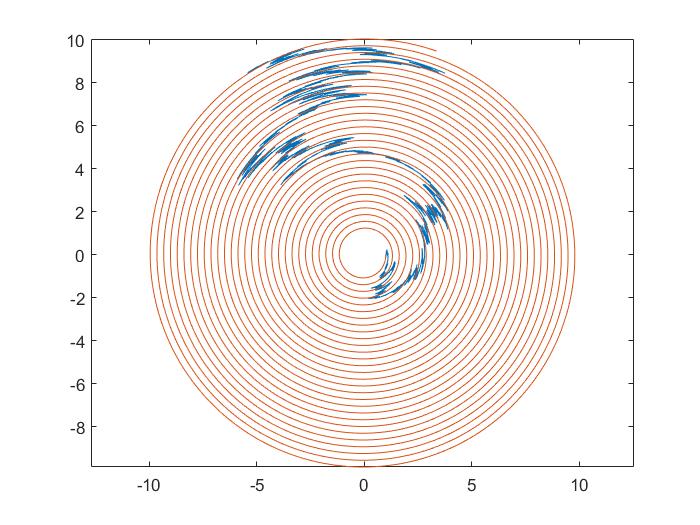}
			\caption{Comparing one sample path of $X$ (blue) and the MCF at the same scale (red).}
			\label{fig:circlemcf}
		\end{figure}
	
		\begin{figure}[h]
			\centering
			\includegraphics[width=0.7\linewidth]{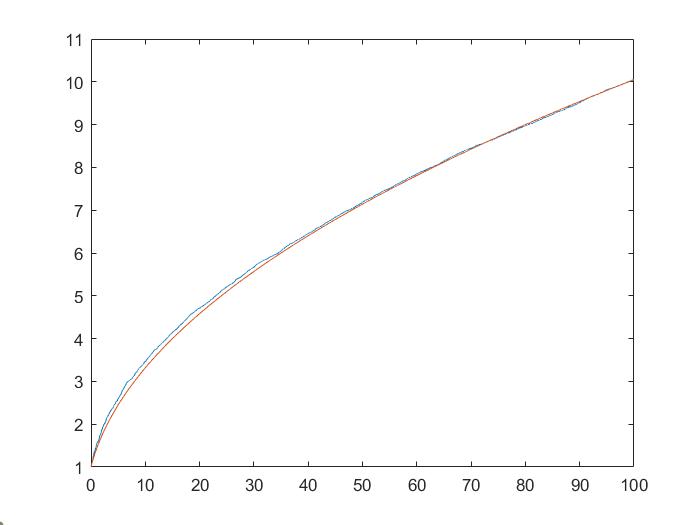}
			\caption{Comparing the theoretical rate (red) of the mean curvature drift and that from the simulation (blue) in \ref{circle example}.}
			\label{fig:mcfrate}
		\end{figure}
		
	\end{example}

		\begin{example}
			Complex projective space with Fubini-Study metric: the Fubini-Study metric on the projective space $\mathbb{CP}^n$ is defined as the quotient metric from the round sphere, i.e. $S^{2n+1} $ embedded in $\mathbb{C}^{n+1}$ with the Euclidean structure as $\mathbb{R}^{2n+2}$, by the unitary group $U(1)$ homeomorphic to the $1$-circle $S^1$. 
			
			It is easy to verify $U(1)$ has totally geodesic orbits, and hence the Laplacian commutes with the quotient. Brownian motion under the Fubini-Study metric is therefore attained by projecting Brownian motion on $S^{2n+1}$, which satisfies
			\begin{equation}  
					dZ = (I-Z^*Z)\strat W = (I-Z^*Z)dW + \frac{1}{2}dZdW = (I-Z^*Z)(dW + \frac{1}{2}dt), \; |Z_0| = 1,
				\end{equation} in homogeneous coordinates.
		\end{example}
	
	\subsection{Brownian motion on orthogonal groups and their quotients}
		Real Grassmannians, flag manifolds, and Stiefel manifolds can be presented as quotients of the orthogonal group $O(n)$, and some most common metrics on them are the quotient metric from $O(n)$. The standard metric on $O(n)$, denoted $g$ for this subsection, can be described either as the bi-invariant from the Killing form or equivalently (up to a scalar) the submanifold metric induced from the usual embedding $O(n) \hookrightarrow \mathbb{R}^{n\times n}$. Let's use the latter as the definition, which can be conveniently computed as the Frobenius norm
		\begin{equation}
			g(B,C) = \tr BC^*.
		\end{equation}
		To write down the invariant Brownian motion, we see that the Lie algebra $\mathfrak{so}(n)$, represented as skew-symmetric matrices, has orthonormal basis 
		\[\{A_{ij} = \frac{E_{ij}-E_{ji}}{\sqrt{2}}  | 1 \leq i < j \leq n\}.\]
		
		Brownian motion, $Q$, on $O(n)$ satisfies 
		\begin{equation}\label{O(n) BM}
			d Q  = Q\strat A, 
		\end{equation} where $A =\frac{1}{\sqrt{2}} \left[W^{ij}\right]_{i,j}$ with $W^{ij} = - W^{ji}$ and $\{W^{ij}| 1\leq i < j \leq n\}$ are independent standard Wiener processes.

		Consider subgroups of the form 
		\begin{equation}
			K = \begin{bmatrix}
				I_{n_0}&    &    &  \\
				& O(n_1)	&   & \\
				&   &  \ddots & \\
				&&&           O(n_r)		
			\end{bmatrix}
		\end{equation} with $n_0 + n_1 + \cdots + n_r = n$ and $n_i\geq 0$. Right multiplication of $K$ on $O(n)$ induces the adjoint map on $\mathfrak{so}(n)$ and is an isometry for $g$ : at $U \in O(n)$, for $V \in K$ and $B,C \in \mathfrak{so}(n)$, 
		\begin{equation}
			\tr (UV(V^*BV))(UV(V^*CV))^* = \tr BC^* = g(UB, UC).
		\end{equation} Hence the quotients by $K$ have totally geodesic orbits. 
		
		\subsubsection{Stiefel manifold}
		
		The compact real Stiefel manifold $Stief(k,n)$ parametrizing all orthonormal frames of $k$-dimensional subspaces of $\mathbb{R}^n$ is 
		\begin{equation}
			\{UI_{k,n}| U \in O(n)\} \cong O(n)/O(n-k),
		\end{equation} where \[I_{k,n} = \begin{bmatrix}
			I_k & \\
			& 0
		\end{bmatrix}.\]
		
		BM on $Stief(k,n)$, denoted S, hence satisfies 
		\begin{equation}
			dS = dQ I_{k,n} = Q \strat A I_{k,n}.
		\end{equation}

		\subsubsection{Grassmannian}
		
		The real Grassmannian $\mathbb{G}r(k,n)$ parametrizing all $k$-dimensional subspaces of $\mathbb{R}^n$  is 
		\begin{equation}
			\mathbb{G}r(k,n) = O(n)/O(k)\times O(n-k).
		\end{equation} It is often practically and conveniently considered as projection matrix, yielding the mapping 
		\begin{equation}
			\phi: U \mapsto U I_{k,n} U^*,
		\end{equation}
		
		Since the metric $g$ is bi-invariant, the fibres of the submersion $\phi: O(n) \twoheadrightarrow \mathbb{G}r(k,n)$, i.e. the cosets of the quotient, have constant volume and therefore are totally geodesic. Consequently, $\phi$ projects Brownian motion on $O(n)$ to Brownian motion on $\mathbb{G}r(k,n)$. If $Q$ is a process satisfying \ref{O(n) BM} , then using the projection matrix presentation, writing $P : = QI_{k,n} Q^*$, Brownian motion on $\mathbb{G}r(k,n)$ satisfies
		\begin{equation}
			dP = d QI_{k,n} Q^* = (\strat Q)I_{k,n}Q^* + QI_{k,n}\strat Q^* . 
		\end{equation}
		But \begin{equation}
			\strat Q = Q \strat A = Q dA + \frac{1}{2}dQdA = QdA - nQdt
		\end{equation}
		since
		\begin{equation}
			dQdA = QdAdA = -2nQdt. 
		\end{equation}
		
		So now, 
		\begin{align}
			dP & = QdAI_{k,n}Q^* - nQI_{k,n}Q^*dt + QI_{k,n}(dA^*)Q^* - nQI_{k,n}Q^*dt \\
			& = -2nPdt + Q(dAI_{k,n} - I_{k,n}dA)Q^*\\
			& = -2nPdt + Q\sum_{i\leq k < j}A_{ij}dW^{ij}Q^* \\
			& = Q(\sum_{i\leq k < j} A_{ij}dW^{ij} - 2nI_{k,n}dt)Q^*. 
		\end{align}
		 
		\begin{remark}
			In \cite{baudoin2020brownian}, they described BM on complex Grassmannians in inhomogeneous coordinates and are able to derive autonomous equation for the process along with detailed eigenvalue processes:   
			
			Write the invariant BM on $U(n)$ defined as $dU = U\strat dA$ (similar to $Q$ above)  in blocks as
			\begin{equation}
				U = \begin{bmatrix}
					X_{(n-k)\times k} & Y_{(n-k)\times (n-k)}\\
					Z_{k\times k} & C_{k \times (n-k)}
				\end{bmatrix},
			\end{equation} with the sizes indicated in the subscripts. In the homogeneous coordinates we use, 
			\begin{equation}
				P := UI_{k,n}U^* = \begin{bmatrix}
					X\\
					Z
				\end{bmatrix}
				\begin{bmatrix}
					X^* & Z^*
				\end{bmatrix}
				= \begin{bmatrix}
					XX^* & XZ^* \\
					ZX^* & ZZ^*
				\end{bmatrix}
			\end{equation} represents the BM on the complex Grassmannian $Gr_{\mathbb{C}}(k,n)$ as a projection matrix.
		
			In \cite{baudoin2020brownian}, they consider the processes $G: = XZ^{-1}$ and $S = GG^*$. One can prove that the stopping time $\tau : = \inf \{t|\det Z_t = 0\} = +\infty$ a.s. $S$ satisfies the autonomous equation 
			\begin{equation}
				dS = \sqrt{I_k + S}dB^*\sqrt{I_k + S}\sqrt{S} + \sqrt{S}\sqrt{I_k + S}dB\sqrt{I_k + S} + 2(n-k -\tr(S))(I_k + S)dt,
			\end{equation} where $dB$ is the standard $k\times k$ complex BM, and the eigenvalues of $S$ satisfy
			\begin{equation}
				d\lambda_i = 2(1+\lambda_i)\sqrt{\lambda_i} d\beta_i +2(1+\lambda_i)(n-2k+1-(2k-3)\lambda_i+ 2\lambda_i(1+\lambda_i)\sum_{j\neq i}\frac{1}{\lambda_j -\lambda_i})dt,
			\end{equation} where $\beta_i$'s are independent real BM.
		\end{remark}
		
		\subsubsection{Flag manifolds}
		
		More general than Grassmannians, a real flag manifold $Flag(n_1, \cdots, n_r)$ parametrizes all \textit{flags}  of linear subspaces $V_1 < V_2 < \cdots < V_r = \mathbb{R}^n$ with $\dim_{\mathbb{R}} V_i = n_1 + \cdots + n_i$. Equivalently, a flag is an orthogonal decomposition $\mathbb{R}^n  =  W_1 \oplus \cdots \oplus W_r$ with $V_i = W_1 + \cdots + W_i$. Hence we can represent 
		\begin{equation}
			Flag(n_1, \cdots, n_r) = O(n)/(O(n_1)\times \cdots \times O(n_r)).
		\end{equation}
		
		To present the flag manifold in matrix form, we may consider the embedding
		\begin{equation}
			Flag(n_1, \cdots, n_r) \hookrightarrow Gr(n_1,n) \times Gr(n_r,n).
		\end{equation} That is, for $Q \in O(n)$ with columns $q_1, \cdots, q_n$, let the matrix 
		\begin{equation}
			Q_k = \begin{bmatrix}
				q_{n_{k-1}+1} & \cdots & q_k
			\end{bmatrix}, 
		\end{equation}  and the image of $Q$ in $Flag(n_1, \cdots, n_r)$ is presented as
		\begin{equation}
			(Q_1^*Q_1, \cdots, Q_r^*Q_r).
		\end{equation}
		
		\subsubsection{Normal metric of isospectral orbits}
		% explain why this section
		Let $\Lambda = \diag(\lambda_1, \cdots, \lambda_n)$. The \textit{adjoint orbit} of $\Lambda$ by $O(n)$ is
		\begin{equation}
			Ad_{O(n)} \Lambda = \{U\Lambda U^*|U \in O(n)\}.
		\end{equation} The \textit{normal metric} on $Ad_{O(n)} $ is defined as the quotient metric from $(O(n), g)$. Indeed, WLOG, assuming $\lambda_1 = \cdots \lambda_{n_1} > \lambda_{n_1 +1} = \cdots \lambda_{n_1+ n_2} > \cdots > \lambda_{n_1+\cdots+n_{r-1} + 1 } = \cdots = \lambda_n$, we see 
		\begin{equation}
			Ad_{O(n)} \Lambda  \cong Flag(n_1, \cdots, n_r).
		\end{equation}
		
		%	\subsubsection{Bures-Wasserstein metric restricted to isospectral orbits}

	\subsection{Metric lift of reductive homogeneous spaces}
	
	It is not completely obvious that the $G$-invariant Riemannian structures on some common homogeneous spaces $G/K$ come from Riemannian submersions--a submersion $\phi: M\to (N,h)$ can only be Riemannian if for each $a\in N$, one can find decompositions $T_m M = Ver_m M \oplus Hor_m M$ for each $m \in \phi^{-1}(a)$ and isometries between all $Hor_m M$ along the fibre. 
	
	We list a series of facts stating that it is always the case, and the metrics come from the $G$-invariant metric on $G$. In this way, understanding the BM on $G/K$ is more or less reduced to that on $G$. In fact, the two examples we shall see in the following sections are \textit{Riemannian symmetric spaces}, where we have canonical metric lifts. The material follows mostly from \cite{gallier2020differential}.

	 \begin{definition}
		%reductive space
		Let $G$ be a Lie group and $K$ a closed subgroup. $Lie(G)  = \mathfrak{g}$ and $Lie(k) = \mathfrak{k}$. The homogeneous space $G/k$ is \textit{reductive} if $\mathfrak{g}$ of $G$ is a direct sum \[ \mathfrak{g} = \mathfrak{h} \oplus \mathfrak{m}\] as vector spaces such that $\mathfrak{m}$ is stable under the adjoint representation of $K$, i.e. \[ Ad_k \mathfrak{m} \subset \mathfrak{m}\] for any $k \in K$.
	\end{definition}
	
	\begin{proposition}\label{reductive lift}
		
		Let the space of right cosets, $G/K$, be a reductive homogeneous space as described in the definition above. 
		\begin{enumerate}
			\item If $K$ is compact, then there exists a left $G$-invariant metric on $G/K$. Furthermore, if the adjoint representation of $H$ on $\mathfrak{m}$ is irreducible, then the metric is unique up to scalars. 
			
			\item If $\mathfrak{m}$ has an $Ad(K)$-invariant inner product $\langle\cdot,\cdot\rangle_\mathfrak{m}$, then is quotient map $G \to G/K$ is a Riemannian submersion where $G$ is endowed with any left $G$-invariant metric $\langle\cdot, \cdot\rangle_\mathfrak{g}$ extending $\langle\cdot,\cdot\rangle_\mathfrak{m}$ such that $\mathfrak{k}$ is orthogonal to $\mathfrak{m}$, where $\mathfrak{m}$ maps isometrically to $T_K G/K$.
		\end{enumerate}
	\end{proposition}
	
	A class of homogenous spaces have a natural reductive decomposition, namely symmetric spaces, introduced by \'Elie Cartan.
	
	\begin{definition}
		Let $G$ be a Lie group. An \textit{involution} or \textit{involutive automorphism} $\sigma$ is an automorphism such that \[\sigma^2 = Id, \; \sigma \neq Id.\] 
		
		We denote $G^\sigma = \{a\in G| \sigma(a) = a\}$ and $G^\sigma_e$ the identity component. It is easy to check $G^\sigma$ is a subgroup. 
	\end{definition}

	\begin{definition}
		Let $G$ be a Lie group and $K$ a closed subgroup. A Cartan involution  $\sigma$ for the homogeneous space $G/K$ is an involutive automorphism of $G$ such that $G^\sigma_e \subset K \subset G^\sigma$.
	\end{definition}

	A Cartan involution gives a natural choice for a reductive decomposition for $G/K$. Note the differential map $\sigma_*$ has eigenvalues $\pm 1$ on $Lie(G) =:\mathfrak{g}$.
	
	\begin{proposition}
		If $\sigma$ is a Cartan involution for the homogeneous space $G/K$, then \[Lie(G) = \mathfrak{g} = \mathfrak{k} \oplus \mathfrak{m}, \] where $Lie(K) = \mathfrak{k}$ is the eigenspace of $\sigma_*$ corresponding to eigenvalue $1$, and $\mathfrak{m}$ corresponds to eigenvalue $-1$.
	\end{proposition}

	\begin{definition}
		Let $\sigma$ be a Cartan involution for the homogeneous space $G/K$. If $G$ is connected, and $G^\sigma_e$ and $K$ are compact, $G/K$ is called a \textit{symmetric space}. 
	\end{definition}

	\begin{proposition}
		Let $G/K$ be a symmetric space with Cartan involution $\sigma$. Then $G/K$ has $G$-invariant metrics. Moreover, the metrics satisfy conditions in \ref{reductive lift} and thus can be lifted to $G$.
	\end{proposition}
	
	In other words, the existence of a Cartan involution gives a natural way to lift the invariant metric to $G$. In fact the quotients of orthogonal groups introduced in the previous sections are also symmetric spaces. However the importance of the following examples are those whose metrics are not, at least to the author's knowledge, traditionally defined from the point of views of quotient metrics and could provide new perspectives.

	\subsection{Brownian motion on symmetric spaces}
	
	\begin{proposition}\label{symmetric space MCF}
		Let $G \to G/K$ be a Riemannian submersion to a symmetric space as described previously. Then the fibres are totally geodesic. 
	\end{proposition}

	\begin{proof}
		As stated in \ref{reductive lift}, the inner product on $\mathfrak{m}$ is $\Ad_K$-invariant. Hence right $K$-multiplication is an isometry by noticing it corresponds to $\Ad_K$ on horizontal spaces, so the mean curvature is constant along each fibre and MCF commutes with the submersion.
		
		Moreover, since the lifted metric on $G$ is left invariant, all fibres $aK$ have equal volumes, and therefore the volume function has zero gradient, meaning the mean curvature is zero by \ref{log vol grad}. 
	\end{proof}
	
 		\subsubsection{Poincar\'e's upper half plane}
		Consider the right-quotient from the special linear group to Poincar\'e upper half plane,with $i = \sqrt{-1}$, \begin{equation}\label{Poincare plane}
			q:	SL(2,\mathbb{R}) \to 	SL(2,\mathbb{R}) /SO(2)\cong \mathbb{H}_+ = \{x + iy| x\in \mathbb{R}, \; y>0\},
		\end{equation} via the action
		
		\begin{equation}
			\begin{bmatrix}
				a & b\\
				c & d
			\end{bmatrix} \cdot i
		=  \frac{ai + b}{ci + d},
		\end{equation} which is an isometry on $\mathbb{H}_+$.
	%Importance: 
		
		The upper half plane $\mathbb{H}_+$ equipped with the Poincar\'e metric is a classical object.
		
		\begin{equation}
			ds^2 = \frac{dx^2 + dy^2}{y^2}
		\end{equation} with the usual $(x,y)$ coordinate. 
		
		Taking $\sigma(M) = M^{-T}$, it is easy to check $\mathbb{H}_+$ is a symmetric space with 
		
		\begin{equation}
		\mathfrak{k} : =	\mathfrak{so}(2,\mathbb{R}) = \langle Z = \frac{1}{2} \begin{bmatrix}
								0 & -1 \\
								1 & 0
							\end{bmatrix}\rangle
		\end{equation}
	
		and \begin{equation}
			\mathfrak{m} : = \langle 
						X = \frac{1}{2}\begin{bmatrix}
								0 & 1 \\
								1 & 0
							\end{bmatrix}, \;			
						Y = \frac{1}{2}\begin{bmatrix}
								1 & 0 \\
								0 & -1
							\end{bmatrix}\rangle.
		\end{equation}
	
		One can check that $q_*: X \mapsto \partial_x \in T_i \mathbb{H}_+$ and $q_*: Y \mapsto \partial_y \in T_i \mathbb{H}_+$, so we can pick any metric on $\mathfrak{sl}(2)$ such that $\mathfrak{m} \perp \mathfrak{so}(2)$ and $\mathfrak{m}$ is isometric to $T_i \mathbb{H}_+$.
		
		A BM on $SL(2)$ under this metric satisfies
		\begin{equation}
			dM = M\strat W
		\end{equation} for 
	
		\begin{equation}
			W = \frac{1}{2}\begin{bmatrix}
				W^Y & W^X-W^Z\\
				W^X + W^Z & -W^Y 
			\end{bmatrix},
		\end{equation} where $W^X$, $W^Y$, $W^Z$ are independent standard Wiener processes. A BM on $SL(2)$, $M = \begin{bmatrix}
		a & b\\
		c & d
	\end{bmatrix}$, projects directly to a BM on $\mathbb{H}_+$, $M\cdot i = \frac{ai+b}{ci +d}$.

	\subsubsection{The Cartan-Hadamard geometry on $\pdef$}\label{C-H BM}

		The \textit{Cartan-Hadamard} or \textit{trace metric} on $\pdef \cong GL(n)/O(n) \cong GL^+(n)/SO(n)$, where $GL^+(n)$ denotes the identity component of $GL(n)$, is \begin{equation}
			\langle V,W\rangle = \tr P^{-1}VP^{-1}W
		\end{equation} at $P \in \pdef$.  One can check it is invariant under the action 
		\begin{equation}
			GL(n)\times\pdef: (G,P) \mapsto GPG^*,
		\end{equation} which is induced from left multiplications of $GL(n)$ on itself through the mapping $\pi: M \mapsto MM^*$.
	
		Again we take the Cartan involution $\sigma: P \mapsto P^{-T}$. So $GL(n)^\sigma = SO(2)$, and $\mathfrak{m} = S(n)$, the space of symmetric matrices. Therefore the C-H metric is in fact the quotient metric from the left invariant metric on $GL(n)$ with the usual inner product on $\mathfrak{gl}(n) \cong \mathbb{R}^{n\times n}$.

		In equation, a BM $G$ on $GL(n)$ here satisfies
		\begin{equation}
			dG = G\strat W = GdW + \frac{1}{2}dGdW = GdW + \frac{1}{2}Gdt\end{equation}
		Write $ P : = GG^*$ for the BM on $\pdef$ under the C-H metric; 
		\begin{equation}
			dP = G(dW+dW^*)G^* + (n+1)Pdt.\end{equation}
		More about this process can be found in \cite{article}.

	\subsection{An example of nonhomogeneous metric: the Bures-Wasserstein geometry}\label{BWBM section}

	Examples of submersions where the fibres are not totally geodesic and how the BM behave can be found in \cite{je1990riemannian}. We present one more example, namely the Bures-Wasserstein geometry. In the next section, we shall consider it as the model example for developing a stochastic algorithm using the drift. 
	 
	\begin{lemma}\label{BW map} (C.f.\cite{massart2020quotient})
		The mapping 
		\begin{equation}
			\begin{array}{cccc}
				\pi_k: & \mathbb{R}^{nk^\times} &  \to  &  \mathcal{S}_+(k,n) \cong \mathbb{R}^{nk^\times}/O(k) \\ 
				& M & \mapsto  & MM^*
			\end{array}
		\end{equation} is a Riemannian submersion for $k\leq n$. 
		
		Here $\mathbb{R}^{nk}$ is the Euclidean space of real $n$-by-$k$ matrices  equipped with the Frobenius norm, $\left<X,Y\right>_{Fr} := \tr XY^*$, and $\mathbb{R}^{nk^\times}$ is the dense open subset of all rank $k$ matrices.  $\mathcal{S}_+(k,n)$ is the set of all $n$-by-$n$ positive semidefinite matrices of rank $k$.
	\end{lemma}
	
	\begin{definition}
		The induced metric on $\mathcal{S}_+(k,n)$ in the lemma above is called the \textit{Bures-Wasserstein (B-W) metric}.
	\end{definition}
	
	The B-W geometry has been studied in several aspects \cite{massart2019curvature}\cite{massart2020quotient}\cite{bhatia2019bures}. The metric has various importance: it induces the Wasserstein distance on Gaussian measures and has deep connection with optimal transport; furthermore, the curvature is positive, which can be useful for optimization algorithms.  Much of the content in this paper is extendable to singular matrices though we focus on the locus of nonsingular matrices $\pdef$ mostly for the time being.

	We are able to write down the orbit mean curvature in detail. 
	\begin{proposition}\label{J formula spectrum}
		If $P = U\Lambda U^*\in \psdef$ for $\Lambda = \diag\left(\lambda_1,\cdots , \lambda_k, 0, \cdots\right)$, $\lambda_i>0$, and $U\in O(n)$, the drift term
		\begin{equation}
			J(P) = U\tilde{\Lambda}U^*, 
		\end{equation} where \begin{equation}
			\tilde{\Lambda} = \diag\left(\sum_{j \neq 1} \frac{\lambda_1}{\lambda_1+\lambda_j} ,\cdots, \sum_{j \neq k} \frac{\lambda_k}{\lambda_k+\lambda_j}, 0, \cdots, 0 \right)
		\end{equation}
	\end{proposition}
	
	\begin{proof}
		Since $O(k)$ acts by isometry, it suffices to compute the mean curvature at $M = U\tilde{L}$ in a fiber, where via SVD, $U$ is orthogonal and $\tilde{L} = \begin{bmatrix}
			L \\
			0
		\end{bmatrix}$ is $n$-by-$k$ for $L = \diag(l_1, \cdots, l_k)$. 
		
		First, we need to find an orthonormal basis (o.n.b.) for $T_MMO(k) = M\mathfrak{so}(k)$: let $\{E_{ij}\}$ denote the standard basis with the $(i,j)$-th entry being one and zeros otherwise, and let $A_{ij} = E_{ij} - E_{ji}$. Now for $A,B \in \mathfrak{so}(k)$ (i.e. skew-symmetric), \begin{equation}\label{inner product simplification}
			<MA,MB>_{T_MMO(n)} = \tr ULAB^*L^*U^* = \tr LAB^*L.
		\end{equation} It is easy to verify that \begin{equation} \label{onb}
			M\{\tilde{A}_{ij} : = \frac{A_{ij}}{\sqrt{l_i^2 + l_j^2}}| i< j\}
		\end{equation} is an o.n.b. for $T_M MO(k)$.
		
		Each $\tilde{A}_{ij}$ has trajectory \begin{equation}
			c_{ij}(t) = Me^{t\tilde{A}_{ij}}
		\end{equation} on $MO(n)$. Hence the second fundamental form 
		\begin{equation}
			\sff(M\tilde{A}_{ij}, M\tilde{A}_{ij}) = (D_{M\tilde{A}_{ij}} M\tilde{A}_{ij})^\perp = (c_{ij}^{''}(0))^\perp = (-M\frac{E_{ii} + E_{jj}}{l_i^2 + l_j^2})^\perp =-M\frac{E_{ii} + E_{jj}}{l_i^2 + l_j^2}.
		\end{equation} %(Where the tangent vector is extended to a vector field for the covariant derivative. According to \ref{inner product simplification}, $c_{ij}''(0)$ is already perpendicular to the tangent space. )
		
		Therefore the mean curvature is \begin{align}
			\dim O(k) H & = -M(\sum_{i<j}\frac{1}{l_i^2 + l_j^2}(E_{ii} + E_{jj}))
			& = -M\diag(\cdots, \sum_{j \neq i} \frac{1}{l_i^2  + l_j^2} , \cdots)\\ % Check!
			& = -U\diag(\cdots, \sum_{j \neq i} \frac{l_i}{l_i^2  + l_j^2} , \cdots)
		\end{align} and 
		\begin{align}\label{J formula}
			\pi_*(H) & = MH^* + HM^* \\
			& = -2M(\sum_{i<j}\frac{1}{l_i^2 + l_j^2}(E_{ii} + E_{jj}))M^* \\
			&  = -2M\diag(\sum_{j \neq 1} \frac{1}{l_1^2  + l_j^2}, \cdots, \sum_{j \neq i} \frac{1}{l_i^2  + l_j^2} , \cdots)M^*\\
			& = 2 U\diag(\sum_{j \neq 1} \frac{l_1^2}{l_1^2  + l_j^2}, \cdots, \sum_{j \neq i} \frac{l_i^2}{l_i^2  + l_j^2} , \cdots)U^*.
		\end{align}
	\end{proof}
	
	With the knowledge of the well-studied Wishart process, we have the complete eigenvalue process of BM under the B-W geometry:
	
	\begin{corollary}\label{BW BM}
		The Brownian motion on $\psdef$ with the B-W metric satisfies the SDE (until leaving $\mathbb{R}^{nk^\times}$) \begin{multline}
			dW_{BW} = d(W W^*) - J(W_{BW})dt \\
			= (dW)W^* + WdW^* + n I dt - J(W_{BW}) dt,
		\end{multline} where $W$ is a standard $n\times k$-dimensional Wiener process. 
		
		The eigenvalues of $W_{BW}$ satisfy
		\[	d\lambda_i =  2\sqrt{\lambda_i}d\beta_i + (n + \sum_{j\neq i}\frac{\lambda_j(3\lambda_i + \lambda_j)}{\lambda_i^2 - \lambda_j^2})dt,\] for $\lambda_i > 0$, where $\beta_i$'s are independent standard real BM.
	\end{corollary}

	\begin{proof}
		The first statement is a direct application of \ref{BMsubmersion}, noting that 
		\begin{equation}
			dWdW^* = nIdt.
		\end{equation}
	
		$WW^*$ is known as the Wishart process \cite{bru1989diffusions}, whose eigenvalues satisfy
		\begin{equation}
			d\lambda_i = 2\sqrt{\lambda_i}d\beta_i + (n + \sum_{j\neq i}\frac{\lambda_i + \lambda_j}{\lambda_i - \lambda_j})dt,
		\end{equation} where $\beta_i$'s are independent standard real BM.
	
		Combining the result from \ref{J formula}, $W_{BM}$ has eigenvalue processes satisfying
		\begin{equation}
			d\lambda_i = 2\sqrt{\lambda_i}d\beta_i + (n + \sum_{j\neq i}\frac{\lambda_j(3\lambda_i + \lambda_j)}{\lambda_i^2 - \lambda_j^2})dt.
		\end{equation}
	\end{proof}
	
	Now we apply \ref{log vol grad} to compute the mean curvature in matrix entries of the usual representation. 
	\begin{lemma}
		The volume function for the fibres of \[\pi:  \mathbb{R}^{nk^\times} \to \pdef\] at $M \in \pi^{-1}(P)$ is  
		\begin{equation}
			vol(MO(k)) = \sqrt{\det g|_{MO(k)}} \cdot vol(O(k)). 
		\end{equation} 
		
		The determinant in the formula is taken under an orthonormal basis of $\mathfrak{so}(k)$, in which case n entry of $g|_{MO(k)}$ is (using Kronecker $\delta$) \begin{equation}
			g_{ij,ml}  = \delta_{ml}P_{jl} + \delta_{jl}P_{im} - \delta_{il}P_{jm} - \delta_{jm}P_{il}.
		\end{equation} 
	\end{lemma}
	
	\begin{proof}
		Applying \ref{volume of orbit}, we may choose an (orthonormal) basis $\{A_{ij} = \frac{E_{ij}-E_{ji}}{\sqrt{2}}| i<j \leq k\}$ for $\mathfrak{so}(k)$. Then \begin{equation}
			vol(MO(k)) = \sqrt{\det g|_{MO(k)}} \cdot vol(O(k)). 
		\end{equation} 
		
		An entry of $g$ is \begin{align}
			g_{ij,ml} & = \tr MA_{ij}A_{ml}^*M^* \\
			& = \delta_{ml}P'_{jl} + \delta_{jl}P'_{im} - \delta_{il}P'_{jm} - \delta_{jm}P'_{il},
		\end{align} with $P' := M^*M$. Since $vol(MO(k))$ is constant on the orbit, we may choose $M = \sqrt{P}$ so that $P' = P$.
	\end{proof}
	
	\begin{remark}
		We may also justify \ref{J formula spectrum} with the gradient formula:
		roughly speaking, the log volume of an orbit $M\O(n)$ in terms of singular values is $\log vol(MO(n)) \propto \frac{1}{2}\sum_{i\neq j}\log(l_i^2+l_j^2)$. Taking gradient in $l_i$'s, we get $\frac{\partial}{\partial l_i} \log vol(MO(n)) \propto \sum_{i\neq j}\frac{l_i}{l_i^2 + l_j^2}$.
	\end{remark}
	
	We do not have a more explicit formula, but we may compute the gradient under the flat metric on $GL(n)$ and project it through the submersion.
	
	\begin{lemma}
		Let $\psi: M \to N$ be a Riemannian submersion and $f: N \to \mathbb{R}$ a smooth function. Then \[ \grad^N f = \psi_*\grad^M (f\circ\psi).\] 
	\end{lemma}
	
	\begin{proof}
		Since $f\circ \psi$ is constant along each fibre, on $M$, \begin{equation}
			\langle \grad f\circ\psi, V\rangle = V(f\circ\psi) = 0
		\end{equation} for any vertical $V$. Therefore $\grad f\circ\psi$ is horizontal. 
		
		On $N$, for any tangent vector $W$, \begin{equation}
			\langle \grad f , W \rangle = \langle \grad f\circ \psi, \bar{W}\rangle = \langle \psi_*(\grad f\circ \psi), W\rangle,
		\end{equation} where $\bar{W}$ is a horizontal lift of $W$. 
	\end{proof}
	
	\begin{example}
		For $GL(2) \to \mathcal{S}_+(2,2)$, write $P = U\diag(\lambda_1,\lambda_2)U^*$, \begin{equation}\label{n=2 J SVD}
			J(P) = U\diag(\frac{\lambda_1}{\lambda_1+\lambda_2}, \frac{\lambda_2}{\lambda_1+\lambda_2})U^* = \frac{1}{\tr P}P.
		\end{equation}
		
		\begin{equation}
			\det g = g_{12,12} = \tr MA_{12}A^*_{12}M^* = \frac{1}{2}\tr P = \frac{\lambda_1 + \lambda_2}{2}.
		\end{equation}
		
		The mean curvature is \begin{equation}
			H = 	-\nabla \log \sqrt{\frac{P_{11} + P_{22}}{2}} = -\frac{1}{2\tr P}\begin{bmatrix}
				2M_{11} & 2M_{12} \\
				2M_{21} & 2M_{22}
			\end{bmatrix} = -\frac{1}{\tr P } M.
		\end{equation} 
		
		Therefore \begin{equation}
			J(P) = \pi_*(H) = -\frac{\dim O(2)}{2}(MH^* + HM^*) = \frac{1}{\tr P} P,
		\end{equation} consistent with \ref{n=2 J SVD}.
	\end{example}

	\begin{remark}
		From Cor. \ref{BW BM}, we see in matrix entries, $W_{BW}$ has the martingale part defined by $dWW^* + WdW^*$ and finite variation part $nIdt - Jdt)$ Meanwhile for the eigenvalue processes, the Wishart process $WW^*$ make one extra contribution. 
		
		The denominator $\lambda_i^2 -\lambda_j^2$ suggests non-collision property for the eigenvalues. The mean curvature drift $\sum_{j \neq i}\frac{\lambda_i}{\lambda_i + \lambda_j}dt$ itself however has an averaging behaviour--if the initial values satisfy $\lambda_1 = \lambda_2 = \cdots$, then they remain equal. We put this into consideration for possible applications, for example when one needs an algorithm to average the ''shape'' of Gaussian distributions. These are reminiscent of two other processes below.
		
		\textit{Dyson's Brownian motion}: consider the space of $n\times n$ hermitian or real symmetric matrices as a Euclidean space over $\mathbb{R}$. If $M$ is a BM on this space, then its eigenvalues satisfy the equation
		\begin{equation}
			d\lambda_i = d\beta_i + \sum_{j \neq i} \frac{dt}{\lambda_i - \lambda_j },
		\end{equation} where $\beta_i$'s are standard BMs. The drift term serves as repulsion and ensures that eigenvalues do not collide. 
		
		\textit{Brownian motion of ellipsoids}: (c.f. \cite{rogers_williams_2000}) let $G$ be a standard right invariant BM on $GL(n)$, i.e. \begin{equation}
			\strat G = (\strat W_{\mathfrak{gl}(n)})G, 
		\end{equation} where $W_{\mathfrak{gl}(n)}$ is a standard $ n \times n$ BM. Define $X = G^*G$ and $Y = G^*G$. Let $\lambda_i$ be the eigenvalues of $X$ and $\gamma_i = \frac{1}{2}\log \lambda_i$, then 
		\begin{equation}
			d\gamma_i = dW_i  + \frac{1}{2}\sum_{j \neq i}\coth(\gamma_i - \gamma_j).
		\end{equation} In fact, the eigenvalues never collide in this case as well. 
		
		By \ref{C-H BM}, $Y$ as above is a BM on $\pdef$ under the Cartan-Hadamard geometry, commonly referred to as Dynkin's BM. The eigenvalues of $Y$ do not possess an explicit formula. However the eigenvectors have limits for $t \to \infty$. 
	\end{remark}

\section{An invariant control system using mean curvature drift of group orbits}
In this section, we investigate ideas to utilize the mean curvature drift into a stochastic flow. In short, on our model example of $\pdef$, first we prove that the vertical process, $X$, of Brownian motion $W$ projects through $\pi$ a deterministic process $P$ that is exactly the mean curvature drift whose derivative is $J(P)$. Define a new metric on $\pdef$, \[\tr_R VW^*: = \left<V,W\right>_R : = \tr RVW^*,\]  for $R\in \pdef$. We may still consider the Brownian motion $W_R$ and orthogonal projection under the new metric to obtain modified drift term $J^R(P)dt$. By varying $R$ along the path, we have a control system. We shall see that it is in fact a Lie-theoretical \textit{invariant control system}, meaning the directions accessible for the evolution is a subset of (left) invariant vector fields. 

Invariant control systems already existed in the literature. For introduction, see \cite{sachkov2000controllability}. We shall only introduce the version tailored for our case for the time being. 

Typical control theory asks whether the system is \textit{controllable}, i.e. whether the \textit{attainable}, or \textit{reachable} set, the set of points in the space reachable via the controls from a point $P_0$, is the whole space. We shall see this is not the case for our system, in which we always increase the Loewner partial order ($A>B$ if $A-B$ is positive definite) as time increases. One other nuance is that the condition on the control could affect the controllability. Originally we considered $R_t^{-1}$ to be Lipschitz. A natural choice from the point of view of invariant system is \textit{regulated} functions, i.e. limits of piecewise constant functions.  The conclusions for the control system so far is organized in \ref{control equivalence} and \ref{main thm}, which we begin proving.

	\subsection{Construction}\label{construction}
	
		In order to create freedom for the flow, we must alter some factor in the equation.  
	
	Since the mean curvature is a consequence from the metric, we consider a ``change of metric'' on $GL(n)$. Define the inner product on $T_M GL(n)$: \begin{equation}
		\tr_R VW^*: = \left<V,W\right>_R : = \tr RVW^*,
	\end{equation} where $R>0$ is a positive matrix in $\pdef$. A new metric on $GL(n)$ induces a new horizontal space w.r.t. the submersion $\pi$ and hence a different mean curvature denoted $H^R$. We denote the new drift on $\pdef$ by $J^R$. Choosing different $R$, we have a control system 
	
	\begin{equation}
		\begin{cases}
			dX = \prv_R dW_R \\
			dP := d\pi(X) = J^R(P),
		\end{cases}
	\end{equation} where $W_R$ is Brownian motion under $\tr_R$, and $\prv_R$ and $J^R$ are the orthogonal projection to vertical directions and the drift term defined as earlier, in this new metric. 
	
	\begin{remark}
		Is the metric change a reasonable enough choice to achieve the control? One may choose to project to other directions not parallel to the horizontal space; essentially it means introducing new metric on the total space. One might naturally choose other diffusion processes to replace the Brownian motion; as mentioned in introduction, this is equivalent to metric changes too. 
		
		We may require the metric change to affect $W$ or $\prv$ only instead of both. However, these do not share similar nice description. 
	\end{remark}
	
	It turns out there are more than one reasons this is a reasonable choice of new metrics. First, 
	\begin{lemma}
		The class $\{ \left<,\right>_R | R\in \pdef \}$ consists of all metrics on $GL(n)$ (i) invariant under $O(n)$ right multiplication such that (ii) the metric tensors are constant in the usual $\mathbb{R}^{n\times n}$ coordinates.  In particular, the mean curvature is still constant along the fiber under (a fixed) $\tr_R$.
	\end{lemma} 

	\begin{proof}
		Notations: write $E_{ab}$ for the $n\times n$ matrix whose $(a,b)$-th entry is $1$ and $0$ otherwise. Let $e^a$ be the row vector with $a$-th entry being $1$ and $0$ elsewhere. $Q = [q^{ij}]$ denotes a general orthogonal matrix, of which $q^i$ is the $i$-th row. 
		
		$\{E_{ab}\}$ is the standard basis for $\Mn \cong \mathbb{R}^{n\times n}$.
		
		Let $g$ be a metric satisfying conditions (i) and (ii). (i) means the metric tensor
		\begin{equation}
			g_{ij, kl} = <E_{ij}, E_{kl}>_g = <E_{ij}Q, E_{kl}Q>_g = \sum_{a,b} q^{ja}q^{lb}g_{ia,kb}.
		\end{equation} for any $Q\in O(n)$.
		
		\textbf{Claim} $g_{ij,kl} \equiv 0$ for $j \neq l$ and $g_{ij,kj}$ does not depend on $j$.
		
		Take $q^{j} = e^{\bar{a}}$ and $q^{l} = e^{\bar{b}}$. Then \begin{equation}
			g_{ij,kl} = g_{i\bar{a},k\bar{b}}
		\end{equation} for any pair $a\neq b$. 
		
		Next, for $\bar{a}<\bar{b}$, take the $(i,j)\times (\bar{a}, \bar{b})$-submatrix of $Q$ to be \[\frac{1}{\sqrt{2}}\begin{bmatrix}
			1 & 1\\
			1 & -1
		\end{bmatrix}.\] Now \begin{align}
			g_{ij,kl} & = (\sum_{a\neq b} q^{ja}q^{lb})g_{ij,kl} + \sum_a q^{ja}q^{la}g_{ia,ka}\\
			& = 0 + (\frac{1}{2}g_{i\bar{a}, k\bar{a}} - \frac{1}{2}g_{i\bar{b}, k, \bar{b}}). 
		\end{align} Since $\bar{a}$ and $\bar{b}$ are arbitrary, this again implies $
		g_{ia, ka}$ in constant in $a$ and $g_{ij,kl} = 0$ for $j \neq l$. 
		
		If $g = \tr_R$ for $R = [r_ij]$, then \begin{equation}
			g_{ij,kl} = \tr RE_{ij}E_{kl}^* = \tr RE_{ij}E_{lk} = \delta_{jl} \tr RE_{ik} = \delta_{jl} r_{ki}.
		\end{equation} Hence $\{\tr_R|R>0\}$ is exactly the class of all metrics satisfying (i) and (ii).
		
	\end{proof}
	
	Second, it is equivalent to the pull-back metric induced by the left multiplication of $GL(n)$: for $R = G^*G$: 
	\begin{equation}
		\left<V,W\right>_R = \tr RVW^* = \tr (GV)(GW)^*  = \left<GV,GW\right>_{Fr}.
	\end{equation} 
	
	\begin{lemma}
		With the notations above, we have an isometry 
		\begin{equation}\label{L_G isometry}
			\begin{array}{cccc}
				L_{G^-1}: & (GL(n), \tr_R) & \to & (GL(n), Fr) \\
				& M & \mapsto & G^{-1}M.
			\end{array}
		\end{equation}
	\end{lemma}

	In other words, we permute the Frobenius metric on $GL(n)$ by left multiplication map , $L_G$, $G\in GL(n)$, which are homeomorphisms naturally and commute with the quotient to left cosets. The flow under the new metric  becomes  \begin{equation}
		\dot{P} = J^{R}(P) = -c\phi_*(L_G^*(H(L_G(X_t))) ).
	\end{equation} We vary $R_t$ to have control at all time $t$. Now the freedom we have is the set of all vector fields \begin{equation}
		\{J^R| R\in \pdef\} =: \mathfrak{a}
	\end{equation} which is $G$ invariant on $\pdef$. 
	
	\subsection{Existence}
	Some technical details need to be taken care of.
	
		\begin{lemma}\label{local Lipschitz}
		The projection onto each orbit, $MO(n)$, $\prv$ is locally Lipschitz, and hence the equation for $X$ \ref{X process} exists and is unique locally. 
	\end{lemma}
	
	The proof is based on writing the orthogonal projection in the form of continuous algebraic Lyapunov equations, which has the following solution:
	
	\begin{proposition}
		Given $P, B \in \mathcal{S}(n)$ , if $P$ is stable, i.e. its eigenvalues have negative real parts, then the continuous algebraic Lyapunov equation \[P^TX + XP + B = 0\] has a unique solution for $X\in \Mn$, \[X = \int_{0}^{\infty}{e^{tP}Be^{tP}dt}.\]
	\end{proposition}
	
	\begin{proposition}
		Let $M\in GL(n, \mathbb{R})$ and $F = F_{MM^*} = MO(n)$ be the orbit of $M$ under the $O(n)$ action. Let $P: = MM^*$ and $P':= M^*M$.
		\begin{enumerate}
			\item The vertical component of $ME \in T_MF$ is \begin{equation}
				MK = M \int_0^\infty e^{-tP'}(P'E - E^*P')e^{-tP'}dt
			\end{equation}
			\item Alternatively, a  vector in the form $EM\in T_MF$ has vertical component \begin{equation}
				SM = \int_0^\infty e^{-tP}(PE^* + EP)e^{-tP} dtM. 
			\end{equation}
			\item Let $MA \in T_M MO(n)$, the second fundamental form for the orbit is 
			\begin{equation*}
				\sff(MA, MA) = MA^2 - M \int_0^\infty e^{-tP'}(P'A^2 - A^2P')e^{-tP'}dt.
			\end{equation*}
			
		\end{enumerate}
	\end{proposition}
	
	\begin{proof}
		
		1. Observe at $M$, a vertical vector is of the form $MK$, where $K$ is skew-hermitian. The horizontal space is thus \begin{equation}
			\{M^{-1}S| \tr MKS^*M^{*} = \tr KS^* = 0,\; \forall  K \text{ skew-symmetric}\} = \{M^{-1}S|S \text{ symmetric}\}.
		\end{equation} In other words, we have the decomposition \begin{equation}
			T_M GL(n) = T_M F_P \oplus N_M F_P = M\mathfrak{so}(n)\oplus M^{-1*}\mathcal{S}(n),
		\end{equation} where $\mathcal{S}(n)$ denotes all $n\times n$ symmetric matrices. 	
		
		We need to find skew-symmetric $K$ and symmetric $S$ such that 
		\begin{equation}\label{right hand} ME = MK + M^{-1*}S.\end{equation}
		\begin{equation}\implies M^*ME = M^*MK　+ S. \end{equation}
		Write $P' = M^*M$, since $S$ is symmetric, 
		\begin{equation}P'E - P'K = (P'E - P'K)^* = E^*P' + KP'\end{equation}
		\begin{equation}\implies P'E - E^*P' = P'K + KP',\end{equation}
		where $K$ is the only unknown. Let $  B = P'E - E^*P'$. It becomes in the standard form of the continuous algebraic Lyapunov equation
		\begin{equation}-P'K - KP' + B = 0,\end{equation} with $-P'$ negative definite.
		
		The solution is 
		\begin{equation}\label{tangent} K = \int_0^\infty e^{-tP'}Be^{-tP'}dt = \int_0^\infty e^{-tP'}(P'E - E^*P')e^{-tP'}dt .
		\end{equation} 	
		
		2. Similarly,
		\begin{equation}\label{left hand}
			T_M GL(n) = T_M F_P \oplus N_M F_P = \mathfrak{so}(n)M^{-1*}\oplus \mathcal{S}(n)M.
		\end{equation} 
		
		In other words, for $EM\in T_M\Mn$, 
		\begin{equation} EM = KM^{-*} + SM, \end{equation}  
		with $K$ skew-symmetric and $S$ symmetric.
		\begin{equation} \implies K = EP - SP.\end{equation}
		
		Since $K$ is skew-symmetric,
		\begin{equation}
			\begin{split}
				-K & = P^*E^* -P^*S^* \\
				& = PE^* -PS \\
				& = -EP + SP.
			\end{split}
		\end{equation}
		
		\[	\implies SP + PS = (PE^* + EP) =: B,\]
		or \begin{equation} -SP - PS + B = 0. \end{equation}
		
		Since $-P$ is negative definite, the solution for $S$ is 
		\begin{equation}\label{cotangentcomp}
			S = \int_0^\infty e^{-tP}(PE^* + EP)e^{-tP} dt.
		\end{equation}
		
		3. For $MA \in T_M MO(n)$, with $A$ skew-symmetric, has a trajectory on $F_P$ \[\alpha(t) = Me^{tA}.\]
		
		The second fundamental form $\sff(MA,MA)$ is the normal component of 
		\[\alpha''(0) = MA^2.\]
		
		Using  \ref{tangent},
		\begin{equation} \sff(MA, MA) = MA^2 - M \int_0^\infty e^{-tP'}(P'A^2 - A^2P')e^{-tP'}dt. \end{equation} 
		
	\end{proof}
	
	Writing down these equations is also helpful for using existing Lyapunov solvers for numerical implementations.

	\begin{proof} (of \ref{local Lipschitz})
		
		Since $\prv(M)$ is linear, it suffices to test on $C \in T_N \Mn$ with $\|C\| = 1$ for $N$ in a neighborhood of $M$ that \begin{equation}
			\|\prv(N) C_N - \prv(M) C_M\| < k \|N - M\| 
		\end{equation} for some local constant $k$.
		
		To see this, for $G^{-1}N$ near $M$ and $G^{-1}C \in T_NGL(n)$ with $\|G^{-1}C\|_{R} = \|C\|_{Fr} = 1$, the vertical component is \begin{equation}
			G^{-1}SN = G^{-1}\int_{0}^{\infty}{e^{-tNN^*}(NC^* + CN^*)e^{-tNN^*}dt}N.
		\end{equation}  
		
		We see that the improper integral and its derivative w.r.t. $t$ are uniformly convergent in a compact neighbourhood of $M$, and hence the projection map is continuously differentiable, which implies local Lipschitz continuity. %(Note we want the projection matrix to be Lipschitz so we only need $C$ to have norm $1$.) 
	\end{proof}

	\begin{lemma}
		If $R^{-1}$ is a  Lipschitz function in $P\in \pdef$, then $\prv_R$ is locally Lipschitz.
	\end{lemma}
	
	\begin{proof}
		Taking square root can be chosen to be Lipschitz. Thus by \ref{L_G isometry}, we only need $\prv$ to be locally Lipschitz and reduce the case to the standard metric.
	\end{proof}

	\subsection{Controllability}

	We examine the controllability of the system and prove the main theorem \ref{main thm} in this section.

		By varying metric on $GL(n)$, we have a control system on the homogeneous space $\pdef$:
	\begin{equation}\label{control sys}
		\begin{cases}
			\dot{P_t} & = J^{R_t}(P_t), \\
			R_t &\in \pdef.
		\end{cases}
	\end{equation} 

		\begin{definition}
		Let $G$ be a Lie group. $Lie(G) = \mathfrak{g}$ is identified with left invariant vector fields. A \textit{left invariant control system} on $G$ is of the form 
		\begin{equation}
			\begin{cases}
				\dot{m}_t = m_t V_t \\
				V_t \in \mathfrak{c}|_{m_t} \subset \mathfrak{g}|_{m_t}.
			\end{cases}
		\end{equation} 
		Here $\mathfrak{c}$ is a subset of the set of invariant vector fields $\mathfrak{g}$, and ``$|_{m_t}$'' denotes their restriction to the tangent space $T_{m_t}G$. 
	\end{definition}

	We denote the subset in $T_P\pdef$ of all possible choices of $\dot{P}$ as $\mathfrak{a}(P) = \{J^R(P)|R>0\}$. As shown in the last paragraphs of \ref{construction},
	
	\begin{lemma}
		$\mathfrak{a}(P)$ is the restriction at $T_P\pdef$ of a collection of vector fields  on $\pdef$, denoted as  $\mathfrak{a}$, invariant under the $GL(n)$ action  descended from left $GL(n)$-multiplication on $GL(n)$.
	\end{lemma}
	
	We organise equivalent descriptions of the control system in the following:
	
	\begin{proposition}%\label{control equivalence}
		The following are equivalent descriptions of the same control system for the evolution $P_t \in \pdef \cong GL(n)/O(n)$ in the sense that the controls give the same possible set of $\{P_t\}$ for a fixed time $t$ and an initial condition $P_0$.
		
		\begin{enumerate}
			\item $\dot{P_t} = M_tC_tM_t^*$, where $M_tM_t^* = P$, and the control $C_t$ is any matrix in $\pdef$ such that the  eigenvalues are of the form 
			\begin{equation}
				(\sum_{j \neq 1}\frac{1}{\lambda_1+\lambda_j}, \cdots, \sum_{j \neq k}\frac{1}{\lambda_k+\lambda_j}, \cdots)
			\end{equation} with each $\lambda_j >0$. Here the choice of which $M_t$ does not change the set of possible $\dot{P_t}$. 
			
			\item $\dot{P_t} = J^{R_t}(P_t)$, where the control is $R_t \in \pdef$, and $J^R$ denotes the drift term as previously under metric $\tr_R$.
			
			\item $\dot{P_t} = G_t^{-*}J(G_t^*P_tG_t)G_t^{-1}$, where the control is any $G_t \in GL(n)$.
		\end{enumerate}
	\end{proposition}

	\begin{proof}% (of \ref{control equivalence})
		2 and 3 were already discussed. We only need to derive the formula in 1.
		
		Since $\mathfrak{a}$ is left invariant, written as symmetric matrices, $\mathfrak{a}(P) = \bar{L}_{M*} \mathfrak{a}(I) = M\mathfrak{a}(I) M^*$ for $P = MM^*$, where $\bar{L}_G$ is the left $GL(n)$ action inherited from left multiplication on $GL(n)$.
		
		$\mathfrak{a}(I)$ is the set of vectors $J(P)$ pulled back from each $P=MM^*$ via $\bar{L}_{M*}$. Written as symmetric matrices, using \ref{J formula spectrum}, 
		\begin{multline} \label{a(I) formula}
			\mathfrak{a}(I) = \{M^{-1}J(MM^*)M^{-*}| M \in GL(n)\} \\
			= \{V\diag(\sum_{j \neq 1}\frac{1}{\lambda_1 + \lambda_j}, \cdots, \sum_{j \neq i}\frac{1}{\lambda_i + \lambda_j}, \cdots)V^*|\lambda_i > 0, \; V\in O(n)\}.
		\end{multline}
	\end{proof}
	
	As a corollary, we can now prove 1 in the main theorem.
	
	\begin{proof}(of \ref{main thm} (1))
		From 1 in \ref{control equivalence}, we see $\dot{P_t}$ is always positive definite. 
	\end{proof}
	
	\begin{remark}
		This indicates that we can not have $\dot{P} = J^R(P) = P_\infty - P$ for arbitrary $P_\infty$. In fact, let \begin{equation}
			\alpha :  (\lambda_1, \cdots, \lambda_n) \mapsto (\sum_{j \neq 1}\frac{1}{\lambda_1 + \lambda_j}, \cdots, \sum_{j \neq i}\frac{1}{\lambda_i + \lambda_j}, \cdots);
		\end{equation}
		then the image is contained in the convex cone spanned by $e_i + e_j$, where $\{e_i\}_{i = 1, \cdots, n}$ is the standard basis for $\mathbb{R}^n$. This is observed via writing 
		\begin{equation}
			(\sum_{j \neq 1}\frac{1}{\lambda_1 + \lambda_j}, \cdots, \sum_{j \neq i}\frac{1}{\lambda_i + \lambda_j}, \cdots) = \sum_{i<j} \frac{e_i + e_j}{\lambda_i + \lambda_j}
		\end{equation} % check
	\end{remark}
	
	Furthermore, the image is open in $\mathbb{R}^n$. 
	
	\begin{lemma}
		For $n > 2$, the Jacobian matrix $d\alpha$ is (strictly) negative-definite everywhere in the positive cone. In particular, $Image(\alpha)$ is open in $\mathbb{R}^n$.
	\end{lemma}
	% This is 2 in main thm
	\begin{proof}
		The Jacobian matrix \begin{equation}
			d\alpha = -\begin{bmatrix} 
				\sum_{j\neq 1}\frac{1}{(\lambda_1+ \lambda_j)^2} & \frac{1}{(\lambda_1 + \lambda_2)^2} & \frac{1}{(\lambda_1 + \lambda_3)^2} & \cdots & \frac{1}{(\lambda_1 + \lambda_n)^2} \\
				
				\frac{1}{(\lambda_1 + \lambda_2)^2} & \sum_{j\neq 2}\frac{1}{(\lambda_2 + \lambda_j)^2} & \frac{1}{(\lambda_2 + \lambda_3)^2} & \cdots & \frac{1}{(\lambda_2 + \lambda_n)^2} \\
				
				\vdots & & \ddots & & \vdots \\
				
				\frac{1}{(\lambda_1 + \lambda_n)^2} & \cdots & \cdots & \cdots & 	\sum_{j\neq n}\frac{1}{(\lambda_n + \lambda_j)^2}
				
			\end{bmatrix}
		\end{equation}

		\begin{claim}
			$d\alpha = -\sum_{k = 1}^{n-1} M_kM_k^*$ is a sum of ``squares'', where $rank(M_k) = n-k$.
		\end{claim} 
		
		The is achieved by taking \begin{equation}
			\tilde{M}_k = \begin{bmatrix}
				0 & \frac{1}{\lambda_k + \lambda_{k+1}} &  \frac{1}{\lambda_k + \lambda_{k+2}} & \cdots &  \frac{1}{\lambda_k + \lambda_{k+n}} \\
				
				& \frac{1}{\lambda_k + \lambda_{k+1}} & & & \\
				
				& &  \frac{1}{\lambda_k + \lambda_{k+2}} & & \\
				
				&&& \ddots & \\
				
				&&&&  \frac{1}{\lambda_k + \lambda_{k+n}}
			\end{bmatrix}
		\end{equation} (zero when not indicated)
		
		and then let \begin{equation}
			M_k = \begin{bmatrix}
				0_{k-1\times k-1} & \\
				& \tilde{M}_k
			\end{bmatrix}
		\end{equation}  
		
		We see that \begin{equation}
			M_kM_k^* = \begin{bmatrix} \sum_{j = k+1}^{n} \frac{1}{(\lambda_k + \lambda_j)^2} & \frac{1}{(\lambda_k + \lambda_{k+1})^2} & \cdots & \cdots & \frac{1}{(\lambda_k + \lambda_{n})^2} \\
				
				\frac{1}{(\lambda_k + \lambda_{k+1})^2} &  \frac{1}{(\lambda_k + \lambda_{k+1})^2} & & & \vdots\\
				
				\vdots & &  \frac{1}{(\lambda_k + \lambda_{k+2})^2} & & \vdots \\
				
				\vdots & & & \ddots & \vdots \\
				
				\frac{1}{(\lambda_k + \lambda_{n})^2} & & &  \frac{1}{(\lambda_k + \lambda_{n})^2}
			\end{bmatrix},
		\end{equation} which sums to $-d\alpha$. Indeed, each off-diagonal entry is nonzero in $M_{min\{k,j\}}$. On the diagonal, the $(k,k)$-th entry is \begin{equation}
			\frac{1}{(\lambda_1 + \lambda_k)^2} + \cdots + \frac{1}{(\lambda_{k-1} + \lambda_k)^2} + \sum_{j = k+1}^{n} \frac{1}{(\lambda_k + \lambda_j)^2} + 0  =  \sum_{j\neq k}\frac{1}{(\lambda_k+ \lambda_j)^2}
		\end{equation}  
		
		Now $d\alpha$ is negative semidefinite. If it has a zero eigenvalue, then there exists a left-eigenvector $x$ such that \begin{equation}
			-\sum_{k = 1}^{n-1} xM_kM_k^*x^* = 0, 
		\end{equation} implying \begin{equation}
			(xM_k)(M_k^*x^*) = 0,
		\end{equation} or \begin{equation}
			xM_k =0
		\end{equation} for each $k$.
		Let \begin{equation}
			x = \begin{bmatrix}
				1 & -1 & -1 &\cdots & -1
			\end{bmatrix}
		\end{equation} be a left null vector of $M_1$, which has nullity $1$. But \begin{equation}
			xM_{n-1} = -\frac{1}{\lambda_{n-1} + \lambda_n} \neq 0,
		\end{equation} a contradiction (unless $n=2$). Hence $d\alpha$ is nonsingular. 
	\end{proof}
	
	Lastly, we introduce some concepts from invariant control theories and make an estimate for the explorable space by the system \ref{control sys}.

	\begin{definition}
		A \textit{piecewise constant control} is such that for $0 = t_0 < t_1 < \cdots < T$, $V_t = V_i \in \mathfrak{c}$ for $t \in [t_i, t_{i+1})$ and some $V_i$.   
	\end{definition}
	
	A path with piecewise constant control is of the form 
	\begin{equation}
		m(t) = m(0)\exp t_1 V_0 \exp t_2 V_1 \cdots \exp (t-t_i) V_i, 
	\end{equation} if $t\in [t_i, t_{i+1})$.
	
	\begin{definition}
		A \textit{regulated control}  $V_t \in \mathfrak{c}$ for $t\in [0,T]$ is a pointwise limit of piecewise constant controls.
	\end{definition} 
	
	\begin{definition}
		The piecewise reachable set $\mathcal{R}^{pc}(m_0)$ is the union of the paths starting from $m_0$ with piecewise constant controls within finite time. Similarly, the regulated reachable set $\mathcal{R}^{reg}(m_0)$ is the same with regulated controls. 
	\end{definition}
	
	The control system \ref{control sys} can be lifted to $GL(n)$
	\begin{equation}
		\begin{cases}
			\dot{M} = -c\vec{H}^R \\
			-\vec{H}^R \in \hat{\mathfrak{a}},
		\end{cases}
	\end{equation} where the supscript $R$ means the mean curvature $\vec{H}$ of the fibers under the metric $\tr_R$ and $\hat{\mathfrak{a}} = GL(n)\cdot (-H)$ is the orbit of negative mean curvature field $-H$ under $GL(n)$ left multiplication, which maps to $\mathfrak{a}$ surjectively. 
	
	A path with piecewise control  
	\begin{equation}
		P(t) = M(0)\exp t_1 V_0 \exp t_2 V_1 \cdots \exp (t-t_i) V_i \exp (t-t_i) V_i^* \cdots \exp t_1 V_0^* , 
	\end{equation} if $t\in [t_i, t_{i+1})$ and $P(0) = M(0)M(0)^*$. 
	
	As we have seen in \ref{a(I) formula}, $\hat{\mathfrak{a}}$ consists of matrices whose SVD only need to satisfy conditions on the singular values. Hence $V_i$'s can always be chosen to be normal with same eigenspaces, and their exponentials commute. In this case, the only thing that needs to be specified is the eigenvalue. 
	
	The author has no knowledge of any existing results on how in general eigenvalues are affected by multiplications, which would provide much more information about the reachable set.

	\begin{proof}
		3. As discussed, we only need to estimate what eigenvalues are allowed in the reachable set. 
		
		Since regulated functions are the limits of piecewise constant paths, we may take a boundary point $c(e_i+e_j)$, corresponding to $\lambda_i + \lambda_ j = 1/c$ and $\lambda_k = \infty$ otherwise. Since diagonal matrices commute, the exponential of these commute.  
	\end{proof}

\bibliographystyle{alpha}
\nocite{*}
\bibliography{ReportRef}
\end{document}